\documentclass{paper}
\usepackage{graphicx}
\usepackage{amsmath,amssymb, amsfonts,bm,tikz,url}
\usepackage{xcolor}

\def \fb{\mathbf{f}}
\def \ub{\mathbf{u}}
\def \wb{\mathbf{w}}
\def \vb{\mathbf{v}}

\def \xb{\bm{x}}
\def \yb{\bm{y}}
\def \zbm{\bm{z}}

\def \hb{\mathbf{h}}

\def \zb{\mathbf{z}}

\def \nub{\boldsymbol{\nu}}

\newcommand{\bfa}{\bm{\alpha}}
\newcommand{\calL}{{\mathcal{L}_{\bfa}}}
\newcommand{\calV}{{\mathcal{V}_{\bfa}}}
\newcommand{\calG}{{\mathcal{G}_{\bfa}}}
\newcommand{\calD}{{\mathcal{D}_{\bfa}}}
\newcommand{\calC}{{\mathcal{C}_{\bfa}}}
\newcommand{\calCs}{{\mathcal{C}_{\bfa}^*}}
\newcommand{\calDs}{{\mathcal{D}_{\bfa}^*}}
\newcommand{\calN}{{\mathcal{N}_{\bfa}}}
\newcommand{\calT}{{\mathcal{T}_{\bfa}}}

\newcommand{\mc}{\mathcal}
\newcommand{\ds}{\displaystyle}

\newcommand{\RR}{\mathbb{R}}
\newcommand{\re}{\mathbb{R}}
\newcommand{\la}{\langle}
\newcommand{\ra}{\rangle}
\newcommand{\abs}[1]{\lvert #1 \rvert}

\newcommand{\wholeint}{\int_{\mathbb{R}^3}}
\newcommand{\mbRn}{{\mathbb{R}^n}}
\newcommand{\mbR}{{\mathbb{R}}}
\newcommand{\mbf}[1]{\mathbf{#1}}
\newcommand{\omg}{{\Omega}}

\newtheorem{theorem}{Theorem}[section]
\newtheorem{proposition}{Proposition}[section]
\newtheorem{remark}{Remark}[section]

\title{Helmholtz-Hodge decompositions in the nonlocal framework -- Well-posedness analysis and applications}

\begin{document}
\maketitle

\vspace{-1cm}
\noindent \textsf{M. D'Elia}
\textsf{\textit{Sandia National Laboratories, Albuquerque, NM, {\tt mdelia@sandia.gov}}}\\
\textsf{C. Flores,} \textsf{\textit{California State University Channel Islands, Camarillo, CA, {\tt cynthia.flores@csuci.edu}}}\\
\textsf{X. Li,} \textsf{\textit{University of North Carolina-Charlotte, NC, {\tt xli47@uncc.edu}}}\\
\textsf{P. Radu,} \textsf{\textit{University of Nebraska-Lincoln, NE, {\tt pradu@unl.edu}}}\\
\textsf{Y. Yu,} \textsf{\textit{Lehigh University, Bethlehem, PA, {\tt yuy214@lehigh.edu}}}

\vspace{.5cm}
\noindent
\textsf{\textbf{Abstract.}}

\noindent Nonlocal operators that have appeared in a variety of physical models satisfy identities and enjoy a range of properties similar to their classical counterparts. In this paper we obtain Helmholtz-Hodge type decompositions for two-point vector fields in three components that have zero nonlocal curls, zero nonlocal divergence, and a third component which is (nonlocally) curl-free and divergence-free. The results obtained incorporate different nonlocal boundary conditions, thus being applicable in a variety of settings.

\vspace{.5cm}
\noindent
\textsf{\textbf{Keywords.}} Nonlocal operators, nonlocal calculus, Helmholtz-Hodge decompositions.

\vspace{.5cm}
\noindent
\textsf{\textbf{AMS subject classifications.}} 35R09, 45A05, 45P05, 35J05, 74B99.

%%%%%%%%%%%%%%%%%%%%%%%%%%%%%%%%%%%%%%%%%%%%%%%%%%%%%%%%%%%%%%%%%%%%%%
%%%%%%%%%%%%%%%%%%%%%%%%%%%%%%%%%%%%%%%%%%%%%%%%%%%%%%%%%%%%%%%%%%%%%%
\section{Introduction and motivation}
Important applications in diffusion, elasticity, fracture propagation, image processing, subsurface transport, molecular dynamics have benefited from the introduction of nonlocal models. Phenomena, materials, and behaviors that are discontinuous in nature have been ideal candidates for the introduction of this framework which allows solutions with no smoothness, or even continuity properties. This advantage is counterbalanced by the facts that the theory of nonlocal calculus is still being developed and that the numerical solution of these problems can be prohibitively expensive.

In \cite{nlvc} the authors introduce a nonlocal framework with divergence, gradient, and curl versions of nonlocal operators for which they identify duality relationships via $L^2$ inner product topology. Integration by parts, nonlocal Poincar\'e inequality are some of the PDE based-techniques that play an important role in the study of nonlocal operators. In this work we continue to explore the structure of nonlocal operators and take additional steps towards building a rigorous framework for nonlocal calculus, by obtaining Helmholtz-Hodge type decompositions for nonlocal functions. In particular, we show that functions with two independent arguments (labeled two-point functions) can be decomposed in three components: a curl-free component, a divergence-free component, and a component which is both, curl and divergence free.  All operators involved in this decomposition are nonlocal. Under suitable boundary conditions the components are unique in appropriate spaces.

Helmholtz-Hodge decompositions play an important role in many areas; see \cite{bhatia2013helmholtz} for a review of applications and recent literature. In the classical (differential) framework the Helmholtz decomposition theorem states that every vector field on $\re^3$ can be decomposed into divergence-free and curl-free components. A Helmholtz-Hodge (or simply, Hodge) decomposition usually refers to a three-term decomposition on bounded domains, in which the vector flow will be the unique sum of divergence-free, curl-free, and a harmonic component assuming that the components satisfy appropriate boundary conditions. These three components would correspond, in the framework of continuum mechanics, to a decomposition of a motion into a rotational component represented by the curl of a function (thus, divergence free), an expansion or contraction which is the gradient of a potential (hence, curl free), and a translation component (which is both, divergence and curl free). These components can be determined from the original vector field, on the entire space, or on bounded domains, subject to boundary conditions, or other constraints.

These decompositions have numerous applications that extend to many fields outside continuum mechanics. Indeed, as the decompositions rely on duality properties of the operators, counterparts of Helmholtz-Hodge decompositions have been obtained in other frameworks. For example, in graph theory there are discrete versions of the curl and divergence operators, so decompositions along the kernels of these operators have been obtained \cite{H2}. Thus nonlocal decompositions could be seen as extensions of discrete counterparts, as well as with classical decompositions (Section \ref{classic_conv}). These connections mirror existing results that show that the nonlocal theory of peridynamics is an upscaling of molecular dynamics \cite{gspl}. In many instances, results from the classical theory can be seen as particular cases of the nonlocal theory where interaction kernels are taken as derivatives of distributions \cite{nlvc}; in the same manner fractional differential operators are also instances of nonlocal operators for certain kernel functions \cite{delia_fractional}. Thus one may envision the nonlocal realm as a bridge between classical (differential), discrete (molecular dynamics) and fractional (anomalous diffusion) frameworks. Moreover, a key advantage of the nonlocal calculus is its flexibility in modeling, as the operators have the ability to capture different physical features, depending on the choice of the interaction kernel that may be space (even time) dependent, symmetric, integrable, or singular etc. By obtaining Helmholtz-Hodge decompositions corresponding to different kernels we gain flexibility in applications where one may want to choose operators that depend on material or phenomenon characteristics.

In \cite{Du2019} the authors obtain a nonlocal Helmholtz-Hodge decomposition for vector fields that are {\it one-point} functions using one-point nonlocal gradients and curls, thus aligning with classical results. Their results are obtained using Fourier analysis for functions in a periodic setting, a framework which could be more restrictive in applications. Here, the decompositions obtained hold for two-point functions where the components satisfy a priori boundary conditions, a setup that is needed in applications. To our knowledge, the results of this paper are the first nonlocal Helmholtz-Hodge domain decompositions that hold in three-dimensions with nonhomogeneous boundary conditions. At a first glance, the two-point dependence of functions may seem counter-intuitive or cumbersome, given that most of functions in applications assign an output to one individual point. However, a different perspective could be provided by the fact that there is no ``independent" way to measure a quantity at a point, that everything is read and measured relative to a point of reference. Additionally, once a one-point function is given, one can easily transform it into a two point function following different approaches (see Section \ref{applic} for more details on the implementation of this decomposition).

The significance and contributions of this paper can be summarized as follows:
\begin{itemize}
\item[$\bullet$] The Helmholtz-Hodge type decompositions for {\it two-point} vector fields along kernels of nonlocal operators offer some of the first contributions in this direction for nonlocal theories. The examples provided here provide a starting point into future theoretical studies, and in applications, as nonlocal models gain more traction and visibility in different models.
\item[$\bullet$] The nonlocal decompositions are {\it kernel-dependent}, thus providing a versatile framework for applications where interactions between neighboring points could be adjusted.
\item[$\bullet$] The convergence results of Section \ref{classic_conv} show connections with the classical theory of differential operators, an aspect that has been investigated in many papers \cite{DuMengesha}, \cite{fos-rad}, \cite{FRW}, \cite{foss2019}, \cite{RTY}, \cite{radu_wells2019}.
\item[$\bullet$] The proof of existence for potential fields given in Theorem \ref{thm:HDecomp}  includes also a generalization of well-posedness results for nonlocal problems; see Remark \ref{gen} in Section \ref{sec:DuMengesha-comparison}.
\end{itemize}

The paper is organized as follows. In Section 2 we introduce the nonlocal operators and their framework for which we will obtain the nonlocal versions of the Helmholtz-Hodge decompositions. The main results of the paper are presented in Section 3, under two sets of boundary conditions. In Section 3 we also present an example for a nonlocal decomposition and perform a convergence analysis to the local setting. We conclude the paper with some connections to existing results for well-posedness of nonlocal problems and a subsection on the implementation of decompositions for two-point functions in the nonlocal setting.

%%%%%%%%%%%%%%%%%%%%%%%%%%%%%%%%%%%%%%%%%%%%%%%%%%%%%%%%%%%%%%%%%%%%%%
%%%%%%%%%%%%%%%%%%%%%%%%%%%%%%%%%%%%%%%%%%%%%%%%%%%%%%%%%%%%%%%%%%%%%%
\section{Nonlocal vector calculus}
In this work, due to the presence of curl operators (acting on vectors only), we limit the analysis to domains in $\mbR^n$ with $n=3$. In the analysis of the Helmholtz decomposition we use the nonlocal vector calculus (NLVC) introduced in \cite{nlvc} and \cite{Lehoucq-Gunzburger} and applied to nonlocal diffusion in \cite{sirev}. This theory is the nonlocal counterpart of the classical calculus for differential operators and allows one to study nonlocal diffusion problems in a very similar way as we study partial differential equations (PDEs) thanks to the formulation of nonlocal equations in a variational setting. In this work we do not consider diffusion only, but we utilize additional nonlocal operators, e.g. the nonlocal curl introduced in \cite{nlvc}, to mimic the local Helmholtz decomposition. The basic concepts of the NLVC and the results relevant to this paper are reported below.

The NLVC is based on a new concept of nonlocal fluxes between two (possibly disjoint) domains; the derivation of a nonlocal flux strictly follows the local definition and it is based on a nonlocal Gauss theorem (the interested reader may find the complete analysis in \cite{nlvc}). %
We define the following nonlocal divergence operators acting on zero-th, first, and second order tensors, respectively:
\begin{align}\label{div_def}
\begin{split}
\big( \calD_{,0} \psi\big)(\xb):=&\wholeint \left( \psi(\xb,\yb)+\psi(\yb,\xb)\right)\bfa(\xb,\yb)d\yb,\\
\big( \calD_{,1} \boldsymbol{\psi}\big)(\xb):=&\wholeint \left( \boldsymbol{\psi}(\xb,\yb) +
                 \boldsymbol{\psi}(\yb,\xb)\right)\cdot \mbf{\bfa }(\xb,\yb)d\yb,\\
\big( \calD_{,2} \mbf{\Psi}\big)(\xb):=&\wholeint \left( \mbf{\Psi}(\xb,\yb)+\mbf{\Psi}(\yb,\xb)\right) \mbf{\bfa }(\xb,\yb)d\yb,
\end{split}
\end{align}
where $\psi:\mbRn\times \mbRn\rightarrow \mbR$ is a scalar function,
$\boldsymbol{\psi}:\mbRn\times\mbRn\rightarrow\mbRn$ is a vector function,
$\mbf{\Psi}:\mbRn\times\mbRn\rightarrow\mbRn\times\mbRn$ is a matrix function, and
$\bfa: \mbRn\times \mbRn\rightarrow \mbRn$ is an antisymmetric function, i.e. $\bfa(\xb,\yb)=-\bfa(\yb,\xb)
$.

The corresponding adjoint operators with respect to the $L^2$ duality pairing
\footnote{As an example, for a vector field $\ub$ and a scalar field $v$, the adjoint of a nonlocal operator $\mathcal A$ is defined as
$
(\mathcal{A}\ub,v)_{\re^n}=(\ub,{\mathcal{A}^{\ast}} v)_{\re^n\times \re^n}.
$}
are defined as \cite{nlvc}
\begin{align}\label{grad_def}
\begin{split}
\big( \calD_{,0}^* \vb\big)(\xb,\yb)=&
\left( \vb(\yb)-\vb(\xb)\right)\cdot \mbf{\bfa }(\xb,\yb),\\[1mm]
\big( \calD_{,1} ^* v\big)(\xb,\yb)=&
\left( v(\yb)-v(\xb)\right)\mbf{\bfa }(\xb,\yb),\\[1mm]
\big( \calD_{,2} ^* \mbf{v}\big)(\xb,\yb)=&
\left( \mbf{v}(\yb)-\mbf{v}(\xb)\right)\otimes \mbf{\bfa }(\xb,\yb).
\end{split}
\end{align}
As in the classical calculus, we define nonlocal gradient operators as the negative adjoint of divergence operators acting on $i$th order tensors, i.e.
\[
\calG_{,i}=- \calD_{,i}^*.
\]
Then, as in the local case, we define the nonlocal diffusion operator $\calL$ acting on a scalar function $u$ as the composition of divergence and gradient operators, i.e.
\begin{equation}\label{diff_def}
- \left(\calL u\right) (\xb)=\calD_{,1}\left(\calD_{,1}^* u\right)(\xb)
= 2 \wholeint  \left( u(\xb)-u(\yb)\right) \big(\mbf{\bfa }\cdot \mbf{\bfa }\big) d\yb.
\end{equation}
%%%
\begin{remark}
For simplicity of notation in the sequel we will drop the numerical subscripts on the operators, unless necessary for clarity, and keep only the reference to the kernel.
\end{remark}

Given a vector field $\ub$, we also define the nonlocal curl operator and its corresponding adjoint as
\begin{equation}\label{curl_adjoint_def}
\begin{array}{l}
\displaystyle\left(\calC \mbf{u}\right)(\xb):=\wholeint \mbf{\bfa } (\xb,\yb)\times
             \left(\mbf{u}(\xb,\yb)+\mbf{u}(\yb,\xb)\right)d\yb\\[3mm]
\displaystyle\left(\calCs \wb\right)(\xb,\yb):=\;\mbf{\bfa }(\xb,\yb)\times
             \left( \wb(\yb)-\wb(\xb)\right).
\end{array}
\end{equation}
For such operators, by substitution, we have
\begin{equation}\label{curlcurl}
\calC\left(\calCs \wb\right)(\xb)
=-2\wholeint \bfa  \times \big[\left(\wb(\yb)-\wb(\xb)\right)\times\bfa  \big]d\yb
\end{equation}
or, equivalently
\begin{equation}\label{curlcurl_identity}
\calC\left(\calCs \wb\right)(\xb)
=2\wholeint \big((\bfa \otimes \bfa)\left(\wb(\yb)-\wb(\xb)\right) -
  \left(\wb(\yb)-\wb(\xb)\right) |\bfa|^2\big)d\yb ,
\end{equation}
where $\bfa$ is short hand notation of $\bfa(\xb, \yb)$. Furthermore, we have the following result
\begin{equation}\label{curlcurl_operators}
\calC\left(\calCs \wb\right)(\xb)=\calD_{,2}(\calDs_{,2}\wb)-\calD_{,0}(\calDs_{,0}\wb).
\end{equation}
Note that, formally, this is the same expression as in the local calculus, i.e. for a vector field $\bf r$, $\nabla\times(\nabla\times{\bf r}) = \nabla(\nabla\cdot{\bf r}) + \nabla^2{\bf r}$, where the latter represents the vector Laplacian.

\begin{figure}
\centering
\includegraphics[width=0.5\textwidth]{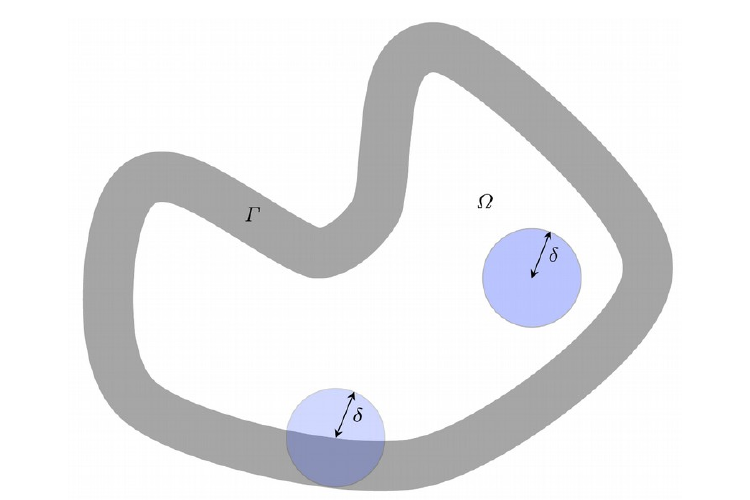}
\caption{The domain $\Omega$, the induced nonlocal boundary $\Gamma$, and balls of horizons $\delta$ centered at points inside the domain $\Omega$ and on the boundary $\partial\Omega$.}
\label{fig:interaction-domain}
\end{figure}

Given an open subset $\omg\subset\mbRn$, the corresponding {\em nonlocal boundary} is defined as
\begin{equation}\label{interaction-domain}
   \Gamma := \{ \yb\in\mbRn\setminus\omg \quad\mbox{such that}\quad \bfa(\xb,\yb)\ne{\bf 0}\quad \mbox{for some $\xb\in\omg$}\}
\end{equation}
so that $\Gamma$ consists of those points outside of $\omg$ that interact with points in $\omg$, see Figure \ref{fig:interaction-domain} for a two-dimensional configuration in which the support of the kernel is a ball of radius $\delta$ as specified in \eqref{eq:kernel-bounded-supp} below.
%%%
%%%
For simplicity, we now focus on the divergence definition {{acting on vectors,}} {{namely $\calD_{,1}$}}, and report important variational results. Corresponding to such a divergence operator we define the action of the nonlocal {\em interaction operator} $\calN(\nub)\colon\mbRn\to\mbR$ on $\nub$ by
\begin{equation}\label{interaction-operator}
\calN\big(\nub\big)(\xb) := -\int_{\omg\cup\Gamma} \big(\nub(\xb,\yb)+\nub(\yb,\xb)\big)\cdot\bfa(\xb,\yb)\,d\yb
\quad\mbox{for $\xb\in\Gamma$},
\end{equation}
and $\calT(\nub)\colon\mbRn\to{\mbR^3}$ by
\begin{equation}\label{interaction-operatorT}
\calT\big(\nub\big)(\xb) := \int_{\omg\cup\Gamma} \big(\nub(\xb,\yb)+\nub(\yb,\xb)\big)\times\bfa(\xb,\yb)\,d\yb,\quad\mbox{for $\xb\in\Gamma$}.
\end{equation}
Note the main difference between local and nonlocal: in the former case the flux out of a domain is given by a boundary integral whereas in the latter case is given by a volume integral.
With $\calD$ and $\calN$ defined as in \eqref{div_def} and \eqref{interaction-operator}, respectively, we have the {\em nonlocal Gauss theorem} \cite{nlvc}
\begin{equation}\label{div-thmn}
\int_\omg \calD(\nub)\,d\xb = \int_{\Gamma} \calN(\nub)\,d\xb.
\end{equation}
Next, let $u(\xb)$ and $v(\xb)$ denote scalar functions; then, the divergence theorem above implies the generalized {\em nonlocal Green's first identity} \cite{nlvc}
\begin{equation}\label{grenns-id}
  \int\limits_\omg v \calD(\calD^\ast u)\,d\xb
- \int\limits_{\omg\cup\Gamma}\int\limits_{\omg\cup\Gamma} \calD^\ast v\cdot \calD^\ast u\,d\yb d\xb=
  \int\limits_\Gamma v \calN(\calD^\ast u)\,d\xb.
\end{equation}
In \cite{nlvc}, one can find further results for the nonlocal divergence operator $\calD$, including a nonlocal Green's second identity, as well as analogous results for nonlocal gradient and curl operators. In addition, in the same paper, further connections are made between the nonlocal operators and the corresponding local operators.

\smallskip
In applications such as mechanics the nonlocal interactions are limited to a finite range, usually much smaller than the size of the domain $\Omega$. Thus, the action of a nonlocal operator at a point $\xb\in\Omega$, the support of the kernel $\bfa$, is limited to a neighborhood surrounding that point; the standard choice of neighborhood consists in Euclidean balls (balls induced by the $\ell^2$ norm) of radius $\delta$, usually referred to as {\it horizon}. More specifically, in this work we consider functions $\bfa:\mbRn\times\mbRn\to\mbRn$ with very mild regularity constraints such that $\bfa(\xb,\yb)\in L^2(\mbRn\times\mbRn)$ and that
\begin{equation}\label{eq:kernel-bounded-supp}
\left\{\begin{array}{ll}
\bfa(\xb,\yb)\geq 0 & \quad\yb\in B_\delta(\xb), \; \xb\in\Omega, \\[2mm]
\bfa(\xb,\yb) = 0   & \quad{\rm elsewhere}.
\end{array}\right.
\end{equation}
As an example, the following $\bfa$ is a very simple example of kernels used to describe the evolution of fractures in peridynamics modeling \cite{silling2010peridynamic}:
\[
\bfa(\xb,\yb)=  { \frac{\yb-\xb}{|\yb-\xb|} \mathcal X(\yb\in B_\delta(\xb))},
\]
where $\mathcal X(\cdot)$ is an indicator function.

For kernels that satisfy \eqref{eq:kernel-bounded-supp}, the associated nonlocal boundary $\Gamma$ (see definition in \eqref{interaction-domain}) could be equivalently defined as
\begin{displaymath}
\Gamma = \{\yb\in\mbRn\setminus\Omega: \exists\;\xb\in\Omega \; {\rm s.t.}\; \|\xb-\yb\|_2\leq \delta\}.
\end{displaymath}

%%%
\noindent \underline{\it Nonlocal boundary conditions.} Due to nonlocality, when solving a nonlocal problem, boundary conditions (i.e. conditions on the solution for $\xb\in\partial\Omega$) do not guarantee the uniqueness of the solution, which can only be achieved by providing conditions on the nonlocal boundary $\Gamma$. As an illustrative example, we consider the following {\it nonlocal Poisson problem} for the scalar function $u$:
\begin{displaymath}
-\calL u = f \quad \xb\in\Omega,
\end{displaymath}
for some $f\in L^2(\Omega)$. Uniqueness of $u$ is guaranteed provided the following conditions are satisfied:
\begin{equation}\label{eq:nBC}
\begin{aligned}
{\rm(D)} & \;\; u = v                & \xb\in\Gamma_D; \\
{\rm(N)} & \;\; -\calN(\calD^*u) = g & \xb\in\Gamma_N,
\end{aligned}
\end{equation}
where $\Gamma_D\cup\Gamma_N=\Gamma$, $\Gamma_D\cap\Gamma_N=\emptyset$, and $v$ and $g$ are some known functions in appropriate functional spaces. The cases $\Gamma_D=\Gamma$ and $\Gamma_N=\Gamma$ are also allowed; however, in the latter case the additional compatibility condition is required: $\int_\Omega f\,d\xb = \int_\Gamma g \,d\xb$; moreover, the uniqueness of solutions is only up to a constant (as in the classical case). We refer to (D) and (N) as Dirichlet and Neumann nonlocal boundary conditions, respectively, due to their resemblance to their local counterparts.

%%%%%%%%%%%%%%%%%%%%%%%%%%%%%%%%%%%%%%%%%%%%%%%%%%%%%%%%%%%%%%%%%%%%%%
%%%%%%%%%%%%%%%%%%%%%%%%%%%%%%%%%%%%%%%%%%%%%%%%%%%%%%%%%%%%%%%%%%%%%%
\section{Nonlocal Helmholtz decomposition}

%%%%%%%%%%%%%%%%%%%%%%%%%%%%%%%%%%%%%%%%%%%%%%%%%%%%%%%%%
\subsection{Existence and uniqueness of the decomposition}
The following proposition is an auxiliary result for the main decomposition theorem.
%%%
\begin{proposition}\label{thm:WP}
Let $\bfa$ be an antisymmetric kernel, $\bfa\in L^2[(\Omega\cup\Gamma)\times (\Omega\cup\Gamma)]$ and suppose $\mc{A}_{\bfa}:L^2(\Omega) \to L^2(\Omega)$ is a positive or negative semidefinite symmetric linear operator such that $\textrm{ker}(\mc{A}_{\bfa})\subset \textrm{ker}({\calCs})$. Then the system
\begin{equation*}%\label{WP_eq1}
\mc{A}_{\bfa} \wb=\mbf{v}
\end{equation*}
is well-posed iff $\mbf{v}=\mc{C}_{\bfa} \mbf{f}$, for some vector field $\mbf{f}\in L^2(\Omega)$.
\end{proposition}
%%%

\medskip\noindent{\it Proof} 
The necessary condition can be simply proved by the fact that
\[
\textrm{Rng}(\mc{A}_{\bfa})^{\perp}
=\textrm{Rng}(\mc{A}^{\ast}_{\bfa})^{\perp}
=\textrm{ker}(\mc{A}_{\bfa})\subset  \textrm{ker}(\calCs)=\textrm{Rng}(\calC)^{\perp}.
\]
where $\textrm{Rng}$ denotes the range of a linear operator, $\perp$ denotes orthocomplementation, and where the second equality is proven in \cite{nlvc}. This implies that $\textrm{Rng}(\mc{A}_{\bfa}) \supset\textrm{Rng}(\calC)$. The sufficient condition is straightforward. $\square$

\medskip
The following theorem is the first nonlocal Helmholtz-Hodge type decomposition that we present in this paper, in which the (gradient) potential function satisfies given nonlocal Dirichlet, respectively, Neumann type boundary conditions.
%%%
\begin{theorem}\label{thm:HDecomp}
For each two-point vector function $\mbf{u}(\xb,\yb)\in L^2(\Omega\times \Omega)$, there exist $\varphi \in L^2(\Omega)$ and $\wb \in L^2(\Omega)$ such that
\begin{equation}\label{HDecomp}
\mbf{u}(\xb,\yb)=\left(\calG \varphi\right) (\xb,\yb) +(\calCs\wb)(\xb,\yb)+ \mbf{h}(\xb, \yb), \quad \xb, \yb \in \Omega
\end{equation}
where $\mbf{h} \in Ker (\calD) \cap Ker(\calC),$ and provided either set of nonlocal boundary conditions are satisfied:
\begin{itemize}
\item[\rm (D)] Dirichlet: $\varphi(\xb)=0$ on $\Gamma$, or
\item[\rm(N)] Neumann: $\calN(\calG(\varphi))=0$ on $\Gamma$ and $\int_\Omega \calD u(\xb)\,d\xb=0$.
\end{itemize}
Uniqueness of $\varphi$ is guaranteed under (D) boundary conditions, while uniqueness up to a constant holds for (N) boundary conditions. Moreover, the decomposition is orthogonal in $L^2(\Omega \times \Omega)$.
\end{theorem}
%%%

\begin{remark}
In the classical Helmholtz decomposition, the counterpart of the term $\mbf{h}$ is a harmonic function, whereas here we have that $\mbf{h} \in Ker (\calD) \cap Ker(\calC)$.
In the nonlocal setting with the operators defined as above, one could not consider $\calL \mbf{h}$, since $\mbf{h}$ is a two-point function. On a related note, the harmonicity of $\mbf{h}$ in the classical decomposition is a consequence of the following results
\[
\text{ curl } \mbf{v} =0\quad \Rightarrow \quad \mbf{ v} =\nabla \psi
\]
and
\[
\text{ div } \mbf{w} =0 \quad \Rightarrow \quad \mbf{ w} =\text{curl } \eta
\]
for appropriate functions $\psi$ and $\eta$ and for domains that are simply connected. However, those results do not hold in the nonlocal case as functions of the form $\calCs\mbf{v}$ do not entirely comprise the null space of the operator $\calD$ (see \cite{nlvc}).

In fact, without additional constraints on $\Omega$ and appropriate boundary conditions, one can easily find nontrivial residual function $\mbf{h}(\xb,\yb)\in Ker (\calD) \cap Ker(\calC)$. For instance,
let
\[
\mbf{h}(\xb,\yb):=\yb-\xb,
\]
then we can confirm from \eqref{div_def} and \eqref{curl_adjoint_def} that $\calD\mbf{h}=\calC\mbf{h}=\bm{0}.$
\end{remark}

%%%%%%%%%%%%PROOF THEOREM 1 %%%%%%%%%%%%%%

\noindent {\it Proof}
We begin by extending the function $\mbf{u}$ by zero
\[
\mbf{u}(\xb,\yb)={\bf 0} \quad \text{ whenever } \quad \xb,\yb \in \re^3 \setminus \Omega.
\]
so that the equality in \eqref{HDecomp} holds on the entire space (similarly, all three terms on the right hand side will be set to zero outside $\Omega\times \Omega$).

We proceed by showing the existence of $\varphi$, followed by $\wb$, while deriving the conditions that the term $\mbf{h}$ must satisfy.

\smallskip{\bf Step 1. Find the (gradient) potential} $\varphi$.
We apply $\calD$ to both sides of \eqref{HDecomp}; we have
\[
\begin{split}
(\calD\mbf{u})(\xb)=&\big(\calD\left(\calG \varphi\right)\big) (\xb) +\big(\calD(\calCs\wb)\big)(\xb)+(\calD\mbf{h})(\xb)\\
=&\big(\calD\left(\calG \varphi\right)\big) (\xb)+(\calD\mbf{h})(\xb)\\
%=&-\big(\calD\left(\calDs \varphi\right)\big) (\xb)+(\calD\mbf{h})(\xb)\\
=&\big(\calL \varphi\big) (\xb)+(\calD\mbf{h})(\xb),
\end{split}
\]
where we used the fact that $\calD(\calCs \wb)=0$. We choose a function $\mbf{h} \in Ker (\calD)$ and let $\varphi$ be the solution to

\begin{equation}\label{phi}
\big(\calL \varphi\big) (\xb)=(\calD\mbf{u})(\xb) \quad \text{ in }\Omega
\end{equation}
under volume constraints (D), or (N). Note that $\varphi$ exists and is unique under condition (D) and it is unique up to a constant for (N).

\smallskip {\bf Step 2. Find} $\wb.$
We apply $\calC$ to both sides of \eqref{HDecomp} and get
\[
\begin{split}
(\calC \mbf{u})(\xb)=&\big(\calC \left(\calG \varphi\right)\big) (\xb) +\big(\calC (\calCs\wb)\big)(\xb)+(\calC \mbf{h})(\xb)\\
=&-\big(\calC\left(\calD^* \varphi\right)\big) (\xb) +\big(\calC (\calCs\wb)\big)(\xb)+(\calC \mbf{h})(\xb)\\
=&\big(\calC (\calCs\wb)\big)(\xb)+(\calC \mbf{h})(\xb)\\
=&-2\wholeint \bfa  \times \big[\big(\wb(\yb)-\wb(\xb)\big)\times\bfa  \big]d\yb +(\calC \mbf{h})(\xb)\\
=&-2\wholeint \big(\wb(\yb)-\wb(\xb)\big) (\bfa \cdot \bfa )d\yb \\
&\qquad +2\wholeint (\bfa \otimes\bfa )\big(\wb(\yb)-\wb(\xb)\big)d\yb +(\calC \mbf{h})(\xb),
\end{split}
\]
where we used the property that $\calC\left(\calD^*\varphi\right)=0$.

We choose a function $\mbf{h} \in Ker(\calC)$ and in order to find $\wb$ we will prove well-posedness for the following equation
\begin{equation}\label{eq:curlWP}
-2\wholeint \hspace*{-.1in}\big(\wb(\yb)-\wb(\xb)\big) |\bfa|^2 d\yb
+2\wholeint (\bfa \otimes\bfa )\big(\wb(\yb)-\wb(\xb)\big)d\yb=(\calC \mbf{u})(\xb),
\end{equation}
where $\mbf{u}(\xb,\yb):\mathbb{R}^n\times \mathbb{R}^n\rightarrow\mathbb{R}^n$ is given. According to \eqref{curlcurl_identity}, \eqref{eq:curlWP} is equivalent to
\[
\calC(\calCs \wb)=\calC\mbf{u}.
\]
The corresponding weak solution $\wb$ then satisfies:
\begin{equation}\label{weakPDE}
\la\calC(\calCs \wb),\mbf{v}\ra_{L^2(\mathbb{R}^n)}=\la\calC\mbf{u},\mbf{v}\ra_{L^2(\mathbb{R}^n)}, \quad \mbf{v}\in \calV,
\end{equation}
where $\calV$ is some proper subspace of $L^2$ that will be defined later on.

Note that, by the definition of adjoint, we have that
\[
\la\calC(\calCs \wb),\mbf{v}\ra_{L^2(\mathbb{R}^n)}=\la \calCs \wb, \calCs\mbf{v}\ra_{L^2(\mathbb{R}^n\times \mathbb{R}^n)},
\]
and
\[
\la(\calC\mbf{u}),\mbf{v}\ra_{L^2 (\mathbb{R}^n)}=\la \mbf{u},\calCs(\mbf{v})\ra_{L^2 (\mathbb{R}^n\times \mathbb{R}^n)}.
\]
Hence, we can first define a bilinear operator $B(\cdot,\cdot)$ as
\[
B(\wb, \wb):=\la(\calCs \wb), (\calCs\wb)\ra_{L^2(\mathbb{R}^n\times \mathbb{R}^n)}.
\]
Note that $B(\cdot,\cdot)$ is not coercive because $K_{\calCs}:={\rm Ker}\left(\calCs\right)$ is not trivial. More specifically, its rank is equal to the one of $K_{\calCs}={\rm span}\{\bfa \}$.

Therefore, we will restrict our function space to be $\calV:=L^2(\mathbb{R}^n)/K_{\calCs}$.
Then on this function space $\calV$, we have the coercivity of $B$:
\[
B(\mathbf{w},\mathbf{w})=\la(\calCs \wb), (\calCs\wb)\ra_{L^2(\mathbb{R}^n\times \mathbb{R}^n)}\ge \|\wb\|_*^2
\]
with $\|\wb\|_*^2:= B(\wb,\wb)$.

Therefore, the weak form \eqref{weakPDE} becomes: find $\wb\in \calV$, such that
\begin{equation}\label{weakPDE_v2}
\la(\calCs \wb), (\calCs\mbf{v})\ra_{L^2(\mathbb{R}^n\times \mathbb{R}^n)}
=\la \mbf{u},(\calCs \mbf{v})\ra_{L^2(\mathbb{R}^n\times \mathbb{R}^n)},
\quad \forall \mbf{v}\in \calV.
\end{equation}

Because for an arbitrary given $\ub\in L^2(\mathbb{R}^n\times \mathbb{R}^n)$, $\la \mbf{u},(\calCs \mbf{v})\ra_{L^2(\mathbb{R}^n\times \mathbb{R}^n)}$ defines a bounded linear functional with respect to any $\mbf{v}\in \calV$, we can apply the Lax-Milgram theorem which guarantees the existence and uniqueness of the solution to \eqref{weakPDE_v2}.
If, instead, we consider $\calV=L^2(\mathbb{R}^n)$, the solution of \eqref{weakPDE}, say $\widetilde{\wb}$, would be unique up to $\zb\in K_{\calCs}$, i.e. $\widetilde{\wb}=\wb+\mbf{z}$. $\square$

\medskip
In the next theorem we prove the existence of a decomposition, provided conditions on the nonlocal normal or tangential component of $\ub$ hold.
\begin{theorem}\label{thm:HDecomp_bc}
For each two-point vector function $\mbf{u}(\xb,\yb)\in L^2(\Omega\times \Omega)$, there exist unique $\varphi \in L^2(\Omega)$ and $\wb \in L^2(\Omega)$ such that
\begin{equation}\label{HDecomp2}
\mbf{u}(\xb,\yb)=\left(\calG \varphi\right) (\xb,\yb) +(\calCs\wb)(\xb,\yb)+ \mbf{h}(\xb, \yb),
\end{equation}
where $\mbf{h} \in Ker (\calD) \cap Ker(\calC),$ $-\calN \calG \varphi=\calN \mathbf{u}$ and $\calT \calCs \wb=\calT \mathbf{u}$. Moreover, the decomposition is orthogonal and unique.
\end{theorem}
%%%

\noindent{\it Proof}
The existence of $\varphi$ and $\wb$ can be proved using similar arguments as in Theorem \ref{thm:HDecomp}. Let $\mathbf{p}=\calG(\varphi(\mathbf{x}))$ and $\mathbf{r}=\calCs (\mathbf{w}(\mathbf{x}))$; by nonlocal integration by parts we have
{\begin{align*}
    {\int_{\omg\cup\Gamma}\int_{\omg\cup\Gamma} (\mathbf{r}+\mathbf{h})\cdot \mathbf{p} \, d\yb d\xb}
& = -\int_{\omg\cup\Gamma}\int_{\omg\cup\Gamma} (\calCs (\mathbf{w})+\mathbf{h})\cdot \calDs(\varphi) d\yb d\xb\\
& = -\int_\Omega \varphi \calD (\calCs (\mathbf{w})+\mathbf{h}) d\xb + \int_\Gamma \varphi \calN (\calCs (\mathbf{w})+\mathbf{h}) d\xb\\
& = 0,
\end{align*}
\begin{align*}
    {\int_{\omg\cup\Gamma}\int_{\omg\cup\Gamma} \mathbf{r}\cdot (\mathbf{p}+\mathbf{h}) \, d\yb d\xb}
& = \int_{\omg\cup\Gamma}\int_{\omg\cup\Gamma} \calCs (\mathbf{w})\cdot
    (-\calDs(\varphi)+\mathbf{h}) d\yb d\xb\\
& = \int_\Omega \mathbf{w}\cdot \calC (-\calDs(\varphi)+\mathbf{h}) d\xb
    -\int_\Gamma \mathbf{w} \cdot \calT (\calCs (\mathbf{w})+\mathbf{h}) d\xb\\
&{=0}.
\end{align*}}
To show the uniqueness of this decomposition, we assume
$$
  \mathbf{u}(\mathbf{x},\mathbf{y})=\calG \varphi_1(\mathbf{x})
+ \calCs \mathbf{w}_1(\mathbf{x}) + \mathbf{h}_1(\mathbf{x},\mathbf{y})%
= \calG \varphi_2(\mathbf{x}) + \calCs \mathbf{w}_2(\mathbf{x})
+ \mathbf{h}_2(\mathbf{x},\mathbf{y}),
$$
then taking a dot product with $(\calG \varphi_1(\mathbf{x})-\calG \varphi_2(\mathbf{x}))$ and integrating over ${\omg\cup\Gamma}$ yield
{\begin{align*}
& \int_{\omg\cup\Gamma}\int_{\omg\cup\Gamma} (\calG \varphi_1(\mathbf{x})
  -\calG \varphi_2(\mathbf{x}))\cdot (\calG \varphi_1(\mathbf{x})-\calG \varphi_2(\mathbf{x})) \\
& +(\calG \varphi_1(\mathbf{x})-\calG \varphi_2(\mathbf{x}))\cdot %
  (\calCs \mathbf{w}_1(\mathbf{x}) + \mathbf{h}_1(\mathbf{x},\mathbf{y})
  -\calCs \mathbf{w}_2(\mathbf{x}) - \mathbf{h}_2(\mathbf{x},\mathbf{y})) d\yb d\xb= 0.
\end{align*}}
By orthogonality and nonlocal integration by parts:
{\begin{align*}
  \int_{\omg\cup\Gamma} &\int_{\omg\cup\Gamma}(\calG \varphi_1(\mathbf{x})-\calG \varphi_2(\mathbf{x}))\\
 & \hspace*{40pt}\cdot   (\calCs \mathbf{w}_1(\mathbf{x}) + \mathbf{h}_1(\mathbf{x},\mathbf{y})-\calCs
    \mathbf{w}_2(\mathbf{x}) - \mathbf{h}_2(\mathbf{x},\mathbf{y})) d\yb d\xb\\
= & -\int_{\omg\cup\Gamma}\int_{\omg\cup\Gamma}\calG \varphi_1(\mathbf{x})
    \cdot(\calCs \mathbf{w}_2(\mathbf{x}) + \mathbf{h}_2(\mathbf{x},\mathbf{y}))d\yb d\xb\\
  &\hspace*{40pt} -\int_{\omg\cup\Gamma}\int_{\omg\cup\Gamma}\calG \varphi_2(\mathbf{x})
    \cdot(\calCs \mathbf{w}_1(\mathbf{x}) + \mathbf{h}_1(\mathbf{x},\mathbf{y}))d\yb d\xb\\
= & \int_\Omega \varphi_1 \calD (\calCs (\mathbf{w}_2)+\mathbf{h}_2) d\xb
    -\int_\Gamma \varphi_1 \calN (\calCs (\mathbf{w}_2)+\mathbf{h}_2) d\xb\\
  &\hspace*{40pt} +\int_\Omega \varphi_2 \calD (\calCs (\mathbf{w}_1)+\mathbf{h}_1) d\xb
    -\int_\Gamma \varphi_2 \calN (\calCs (\mathbf{w}_1)+\mathbf{h}_1) d\xb = 0.
\end{align*}}
Therefore, $\displaystyle\int_{\omg\cup\Gamma}\int_{\omg\cup\Gamma} |\calG \varphi_1(\mathbf{x})-\calG \varphi_2(\mathbf{x})|^2d\yb d\xb= 0$ and {$\calG \varphi_1(\mathbf{x})=\calG \varphi_2(\mathbf{x})$ a.e.}.
Similarly, one can show that {$\calCs \mathbf{w}_1(\mathbf{x})=\calCs \mathbf{w}_2(\mathbf{x})$ a.e.} The above two equalities then yield {$\mathbf{h}_1(\mathbf{x},\mathbf{y})=\mathbf{h}_2(\mathbf{x},\mathbf{y})$ a.e.}. $\square$

%%%%%%%%%%%%%%%%%%%%%%%%%%%%%%%%%%%%%%%%%%%%%%%%%%%%%%%%%%%%%%%%%%%%%%
\subsection{Example of nonlocal Helmholtz decomposition of a field $\bf u$.}\label{sec:example}
In this section we provide an example of decomposition; we focus on functions that live on the plane $\mathcal P=\{\xb\in\Omega\subset\mbR^3: \xb_3=0\}$. The two-point vector function $\mbf{u}: \Omega\times\Omega \rightarrow \mathbb{R}^3$, the one-point scalar function $\varphi:\Omega\rightarrow \mathbb{R}$, and the one-point vector function $\mbf{w}: \Omega\rightarrow \mathbb{R}^2\subset \mathbb{R}^3$ (that is, the range of $\mbf{w}$ could be embedded into $\mathbb{R}^3$).
Hence, according to the definitions of nonlocal operators \eqref{grad_def} and \eqref{curl_adjoint_def} we have
\[
(\calG\varphi)(\xb,\yb): \Omega\times\Omega \rightarrow \mathbb{R}^3\quad
\text{ and }\quad (\calCs \mbf{w}):\Omega\times\Omega\rightarrow
\left(\begin{array}{c} 0\\0\\ \mathbb{R}\end{array}\right).
\]
These test functions of $\varphi$ and $\mbf{w}$ are chosen to be
\begin{equation}\label{test_ex1}
\varphi(\xb):=x_1^2\quad \text{ and }\quad \mbf{w}:=\left(\begin{array}{c} 0\\ x_2^2 \end{array}\right),
\end{equation}
and the nonlocal kernel $\bfa(\xb,\yb)$ is chosen to be
\[
{
\bfa(\xb,\yb)=\frac{1}{\delta^{3/2}}\frac{\yb-\xb}{\left| \yb-\xb\right|}  \mathcal X(\yb\in B_\delta(\xb))
=\frac{1}{\delta^{3/2}|\yb-\xb|}\left(\begin{array}{c}y_1-x_1\\ y_2-x_2\\ 0\end{array}\right) \mathcal X(\yb\in B_\delta(\xb))},
\]
where {{$ B_\delta(\xb)$ denotes a disk on the plane.}}

Consequently, for $\calG\varphi$ we get
\[
(\calG\varphi)(\xb,\yb)
=\frac{-1}{\delta^{3/2}}\frac{\yb-\xb}{\left| \yb-\xb\right|}\left(y_1^2-x_1^2\right)
=\frac{1}{\delta^{3/2}\left| \yb-\xb\right|}
\left(\begin{array}{c}-(x_1^2-y_1^2)(x_1-y_1)\\-(x_1^2-y_1^2)(x_2-y_2)\\ 0\end{array}\right).
\]
For $\calCs \mbf{w}$, it gives
\[
(\calCs\mbf{w})(\xb,\yb)=\frac{1}{\delta^{3/2}|\yb-\xb|}
\left(\begin{array}{c} 0\\ 0\\ (y_1-x_1)(y_2^2-x_2^2) \end{array}\right).
\]
The proposed $\mbf{u}(\xb,\yb)=\calG\varphi+\calCs\mbf{w}$ is equal to
\begin{equation}\label{test_ex2}
\mbf{u}(\xb,\yb)
=\frac{1}{\delta^{3/2}|\yb-\xb|}
\left(\begin{array}{c}-(x_1^2-y_1^2)(x_1-y_1)\\-(x_1^2-y_1^2)(x_2-y_2)\\ (y_1-x_1)(y_2^2-x_2^2) \end{array}\right).
\end{equation}
Now for given $\mbf{u}$ in \eqref{test_ex2}, we will follow the proof of Theorem~\ref{thm:HDecomp} to
solve for the pair $\varphi$ and $\mbf{w}$ defined in \eqref{test_ex1}.

\begin{enumerate}
\item{\textbf{Find $\varphi$}. The corresponding nonlocal equation is,
\[
{{
\begin{cases}
&(\calL\varphi)(\xb)=\ds\frac{2}{\delta^{3}}\ds\int_{B_{\delta}(\xb)}(x_1^2-y_1^2) \,d\yb, \text{ in }\Omega,\\
&\text{volumetric Dirichlet B.C. for}\; \varphi,
\end{cases}
}}
\]
where $d\yb$ is restricted on the plane $\mathbb{R}^2$.
}
\item{\textbf{Find $\mbf{w}$}. The corresponding nonlocal equation is,
\[
\begin{cases}
&\begin{split}
-\frac{2}{\delta^{3}}\ds\int_{B_{\delta}(\xb)}&\left(\mbf{w}(\yb)-\mbf{w}(\xb)\right)d\yb
+\frac{2}{\delta^{3}}\ds\int_{B_{\delta}(\xb)}\frac{\yb-\xb}{|\yb-\xb|}
\otimes\frac{\yb-\xb}{|\yb-\xb|}d\yb\\
% &=\frac{2}{\delta^{3}}\left(\begin{array}{c}\ds\int_{B_{\delta}(\xb)}\frac{1}{|\yb-\xb|^2}(y_1-x_1)(y_2-x_2)(y_2^2-x_2^2)d\yb\\
% \ds\int_{B_{\delta}(\xb)}\frac{-1}{|\yb-\xb|^2}(y_1-x_1)^2(y_2^2-x_2^2)d\yb\\
% 0
% \end{array}
% \right){\text{ in } \mathbb{R}^3} \\
&=\frac{2}{\delta^{3}}\left(\begin{array}{c}\ds\int_{B_{\delta}(\xb)}\frac{1}{|\yb-\xb|^2}(y_1-x_1)(y_2-x_2)(y_2^2-x_2^2)d\yb\\
\ds\int_{B_{\delta}(\xb)}-\frac{1}{|\yb-\xb|^2}(y_1-x_1)^2(y_2^2-x_2^2)d\yb
\end{array}
\right),\text{ in } {\re^2},
\end{split}\\
&\text{volumetric Dirichlet B.C. for }\; \mbf{w},
\end{cases}
\]
where $d\yb$ is again restricted on the plane $\mathbb{R}^2$.}
\end{enumerate}

%%%%%%%%%%%%%%%%%%%%%%%%%%%%%%%%%%%%%%%%%%%%%%%%%%%%%%%%%%%%%%%%%%%%%%
\subsection{Convergence to the local limit}\label{classic_conv}
In this section we study the limit of the operators involved in the decomposition as the extent of nonlocal interactions vanishes. To this end, we consider kernels with support on a ball or radius $\delta$, called horizon or interaction radius, and we study the convergence behavior as $\delta\to 0$. For simplicity and clearness in the exposition, we consider an integrable constant kernel $\gamma = \bfa\cdot\bfa =|\bfa|^2$, such that
\begin{equation}
\bfa(\xb,\yb) = \dfrac{\yb-\xb}{|\yb-\xb|} \mathcal X(\yb\in B_\delta(\xb)).
\end{equation}
We study the limiting behavior of the operators in \eqref{phi} and \eqref{weakPDE}, i.e. $\calL$ and $\calC(\calC^*)$ respectively. While the limit of the nonlocal Laplacian $\calL$ has been widely studied \cite{DuMengesha}, no results have been proved on the limit of $\calC(\calC^*)$; thus, we proceed by analyzing the integral in \eqref{curlcurl_identity}. First, we note that the second term is the vector Laplacian operator, for which we already know that
for some scaling $k_{\delta}$
\begin{equation*}
k_{\delta}\calL \wb(\xb) = \Delta \wb(\xb) + \mathcal O(\delta^2).
\end{equation*}
More specifically, for $\wb=(w_1,w_2,w_3)$, the first component of the vector Laplacian (the other two are obtained in the same way) is given by
\begin{equation*}
-\calL w_1(\xb) = -\dfrac{4\pi}{15}\delta^5(w_{1,x_1x_1}+w_{1,x_2x_2}+w_{1,x_3x_3}) +
                  \textrm{\small HOT},
\end{equation*}
so
\begin{equation*}
\dfrac{15}{4\pi\delta^5 }\calL w_1(\xb) =\Delta w_1(\xb) +
                  \textrm{\small HOT}.
\end{equation*}
Thus, we analyze the other term in \eqref{curlcurl_identity}, i.e. $\ds\int_\mbRn \bfa\otimes\bfa(\wb(\yb)-\wb(\xb))\,d\yb$. We expand the integrand: for any vector $\vb\in\mbR^3$, we have
\begin{equation*}
\bfa\otimes\bfa \vb = \left[
\begin{array}{l}
\alpha_1^2 v_1    + \alpha_1\alpha_2 v_2 + \alpha_1\alpha_3 v_3 \\[1mm]
\alpha_2\alpha_1 v_1 + \alpha_2^2 v_2     + \alpha_2\alpha_3 v_3 \\[1mm]
\alpha_3\alpha_1 v_1 + \alpha_3\alpha_2 v_2 + \alpha^2_3 v_3
\end{array}\right].
\end{equation*}
Due to symmetry, we only study the first component; we have
\begin{equation*}
\begin{array}{l}
\displaystyle2\int\limits_\mbRn \big[\alpha_1^2(w_1(\yb)-w_1(\xb)) +
                  \alpha_1\alpha_2(w_2(\yb)-w_2(\xb)) +
                  \alpha_1\alpha_3(w_3(\yb)-w_3(\xb))\big] \,d\yb  \\[4mm]
\qquad =I + II + III.
\end{array}
\end{equation*}
We analyze each term separately.
\begin{equation*}
\begin{aligned}
I & = 2\int\limits_\mbRn \dfrac{(y_1-x_1)^2}{|\yb-\xb|^2}(w_1(\yb)-w_1(\xb))d\yb \\[2mm]
  & = 2\int\limits_\mbRn \dfrac{h_1^2}{|\hb|^2} (w_1(\xb+\hb)-w_1(\xb))d\hb \\[2mm]
  & = 2\int\limits_\mbRn \dfrac{h_1^2}{|\hb|^2} (\nabla w_1(\xb)\cdot\hb
    + \frac12 \hb^T\nabla^2 w_1(\xb)\hb + \textrm{\small HOT})d\hb.
\end{aligned}
\end{equation*}
It can be shown that, due to symmetry, the first derivative term has no contribution, in fact, the integral is 0. We analyze the term with the Hessian:
\begin{equation*}
\begin{aligned}
2 & \int\limits_\mbRn \dfrac{h_1^2}{|\hb|^2}\;\frac12 \hb^T\nabla^2 w_1(\xb)\hb \,d\hb \\[2mm]
=  & \int\limits_\mbRn \dfrac{h_1^2}{|\hb|^2} (h_1^2 w_{1,x_1x_1} + 2h_1h_2 w_{1,x_1 x_2} + 2h_1h_3 w_{1,x_1 x_3} \\[2mm]
   & \quad + h^2_2 w_{1,x_2 x_2} + 2h_2h_3 w_{1,x_2 x_3} + h_3^2 w_{1,x_3 x_3}) \,d\hb \\[2mm]
=  & \;A+B+C+D+E+F.
\end{aligned}
\end{equation*}
We treat each term separately. It can be shown that $B=C=E=0$; furthermore
\begin{equation*}
A = \dfrac{4\pi}{25} \delta^5 w_{1,x_1 x_1}, \;\;
D = \dfrac{4\pi}{75} \delta^5 w_{1,x_2 x_2}, \;\;
F = \dfrac{4\pi}{75} \delta^5 w_{1,x_3 x_3}.
\end{equation*}
Applying a similar procedure to $II$ and $III$, we have the following:
\begin{equation*}
\begin{aligned}
I&+II+III    \\[2mm]
 &         = \dfrac{4\pi}{25} \delta^5 w_{1,x_1 x_1}
           + \dfrac{4\pi}{75} \delta^5 w_{1,x_2 x_2}
           + \dfrac{4\pi}{75} \delta^5 w_{1,x_3 x_3}
           + \dfrac{8\pi}{75} \delta^5 w_{2,x_1 x_y}
           + \dfrac{8\pi}{75} \delta^5 w_{3,x_1 x_3} \\[2mm]
 &         = \dfrac{4\pi}{75} \delta^5 \Delta w_1 + \dfrac{8\pi}{75} \delta^5 \nabla(\nabla\cdot\wb).
\end{aligned}
\end{equation*}
In summary, the first component of $\calC(\calC^*\wb)(\xb)$ is given by
\begin{equation*}
\begin{aligned}
\big(\calC(\calC^*\wb)\big)_1
& = \dfrac{8\pi}{75} \delta^5 \Big[\big(\nabla(\nabla\cdot\wb)-\Delta w_1\big)-\Delta w_1\Big]
  + \textrm{\small HOT} \\[2mm]
& = \dfrac{8\pi}{75} \delta^5 \Big[\nabla\times(\nabla\times\wb) - \Delta w_1\Big]
  + \textrm{\small HOT}
\end{aligned}
\end{equation*}
Note that this is not consistent with the first component of the local operator, thus we conclude that
\[
\kappa\,\calC(\calC^*\wb)\to \nabla\times(\nabla\times \wb)-\Delta\wb,
\]
where $\kappa=\ds\frac{75}{8\pi \delta^5}$.

In conclusion, the convergence of the scaled nonlocal operator $\kappa\,\calC(\calC^*\wb)$ to its classical counterpart $\nabla\times(\nabla\times \wb)$ holds only on the space of harmonic functions $\wb$.

%%%%%%%%%%%%%%%%%%%%%%%%%%%%%%%%%%%%%%%%%%%%%%%%%%%%%%%%%%%%%%%%%%%%%%
\subsection{Connections to other results in the literature}\label{sec:DuMengesha-comparison}
In \cite{DuMengesha}, the well-posedness of the linearized peridynamics equilibrium system
\begin{equation}\label{eq:linPDWP}
-\int\limits_{B_{\delta}(\xb)\cap\Omega}{\mathbb{C}(\yb-\xb)\left(\wb(\yb)-\wb(\xb)\right)}d\yb =\mbf{b}(\xb),
\end{equation} is established for $\delta>0$,  where $\wb(\xb)$ denotes a displacement field, and where $\mbf{b}(\xb)$ is a given loading force density function and $\mathbb{C}(\bm{\xi})$ is the micromodulus tensor defined by \begin{equation}\label{eq:micromod}
\mathbb{C}(\bm{\xi})=2\frac{\rho(\abs{\bm{\xi}})}{\abs{\bm{\xi}^2}}\bm{\xi}\otimes \bm{\xi} + 2F_0(\abs{\bm{\xi}})\mathbb{I}.
\end{equation} The functions $\rho$ and $F_0$ are given radial functions and their properties determine the well-posedness of \eqref{eq:linPDWP}. We mention the works of \cite{DuMengesha} and others (see for example \cite{Mengesha2010,Du2013c,Du2011b,Emmrich2007b,Emmrich2007a,Lehoucq-Gunzburger,Du2010c}) for well-posedness theory in the case where $F_0\equiv 0,$ although \cite{Silling2000} argues that the condition $F_0\equiv 0$ is too restrictive for equations of motion for bond-based materials. Moreover, equation \eqref{eq:linPDWP} defines a nonlocal boundary value problem with Dirichlet-type volumetric boundary conditions and approaches the Navier equations of elasticity with Poisson ration $1/4$ as $\delta\rightarrow 0$ (see \cite{Emmrich2007b,Du2013c}). Well-posedness is studied in the function space
\vspace*{-25pt}
\begin{multline}\label{eq:funcSpace}
\mathcal S(\Omega):=\left\{ \wb\in L^2(\Omega:\RR^d) :\phantom{\int\limits_{\Omega\times\Omega}} \right.\\[-15pt]
\left.\int\limits_{\Omega\times\Omega} \rho(\abs{\yb-\xb})\left|\frac{\yb-\xb}{\abs{\yb-\xb}}\cdot\left( \wb(\yb)-\wb(\xb)\right)\right|^2  d\yb d\xb < \infty\right\}.
\end{multline}
Whenver $\rho$ is singular, $\mathcal{S}(\Omega)$ is equivalent to a fractional Sobolev space, a proper subspace of
$L^2(\Omega)$ (see \cite{nlvc} and references therein).

We can reformulate \eqref{eq:curlWP} to \eqref{eq:linPDWP} by selecting appropriate $\bfa$. Particularly, we are interested in the prototype kernels $$\bfa = \frac{\yb-\xb}{\abs{\yb-\xb}^{1+\beta}} \mathcal X(\yb\in B_\delta(\xb)),$$ where $\beta=\beta(d)>0$. By letting $\bm{\xi}=\yb-\xb$, the formulations \eqref{eq:curlWP} and \eqref{eq:linPDWP} agree and a simple calculation finds \begin{equation}\label{eq:rho-F0}\rho(\abs{\bm{\xi}})=\frac{1}{\abs{\bm{\xi}}^{\beta}}=F_0(\abs{\bm{\xi}}).\end{equation}

For reference, the simple calculation is included below. Consider the tensor $$\bfa \otimes \bfa - \bfa \cdot \bfa \mathbb{I}.$$ Letting $\bm{\xi}=\yb-\xb$, we have
\begin{align*}
\bfa \otimes \bfa - \bfa \cdot \bfa \mathbb{I} &= \frac{\yb-\xb}{\abs{\yb-\xb}^{1+\beta}} \otimes \frac{\yb-\xb}{\abs{\yb-\xb}^{1+\beta}} - \frac{(\yb-\xb)\cdot (\yb-\xb)}{\abs{\yb-\xb}^{2(1+\beta)}} \mathbb{I} \\
&= \frac{1}{\abs{\bm{\xi}}^{2(1+\beta)}} \bm{\xi}\otimes\bm{\xi} +\frac{\bm{\xi}\cdot \bm{\xi}}{\abs{\bm{\xi}}^{2(1+\beta)}} \mathbb{I}\\
&= \frac{1}{\abs{\bm{\xi}}^{2}} ~ \abs{\bm{\xi}}^{-\beta} \bm{\xi}\otimes\bm{\xi} + \abs{\bm{\xi}}^{-\beta}\mathbb{I}\\
&= \frac{1}{\abs{\bm{\xi}}^{2}} \rho(\abs{\bm{\xi}}) \bm{\xi}\otimes\bm{\xi} + F_0(\abs{\bm{\xi}})\mathbb{I} = \frac{1}{2}\mathbb{C}(\bm{\xi}) .
\end{align*}

\begin{remark}\label{gen}

The identity \eqref{eq:rho-F0} implies that the results in \cite[see Thms 4.2 \& 4.5]{DuMengesha} cannot be used to obtain unique solutions to \eqref{eq:curlWP} for a given $\mbf{u}$. That is, the results of \cite{DuMengesha} apply in two scenarios. The first scenario\footnote{Note that the assumption $|\bm{\xi}|^2\rho(\bm{\xi})\in L^1_{loc}(\RR^d)$ is a simplification of the requirement that $\lim\limits_{\delta\rightarrow 0} \frac{\delta^2}{\ds\int_{B_{\delta}(0)} |\bm{\xi}|^2\rho(\bm{\xi}) d\bm{\xi}} =0. $} requires $F_0\in L^1_{loc}(\RR^d)$ while $|\bm{\xi}|^2\rho(\bm{\xi})\in L^1_{loc}(\RR^d)$. This would require $\beta < d+2,$ while simultaneously having $\beta< d$. In the Sobolev scale, we take $\beta = d+2s,$ for $s\in(0,1)$, hence this scenario is not fruitful. The second requires $F_0$ to have zero mean value and for $\rho\in L^1_{loc}(\RR^d)$ and does not apply to the work presented here. Therefore, the results from Theorems \ref{thm:HDecomp} and \ref{thm:HDecomp_bc} lift previous restrictions mentioned in \cite{Silling2000}.
\end{remark}

%%%%%%%%%%%%%%%%%%%%%%%%%%%%%%%%%%%%%%%%%%%%%%%%%%%%%%%%%%%%%%%%%%%%%%
\subsection{Implementation of results in applications}\label{applic}
In this paper, we produced the Helmholtz decomposition for two-point functions $\ub(\xb,\yb)$. To motivate the applicability of the results, we offer some comments regarding the implementation.

First, the input of a function could be seen as a two-variable argument in which one of the points being fixed, as a point is identified by being referenced to another one (usually, taken to be the origin). Thus, $\fb(\xb)$ is really $\fb(\xb,{\bf 0})$, so our results could be seen as decompositions for functions evaluated at $\xb$, when $\yb$ is the reference point, which would be generalizations of results in the conventional setting where we only track the variable $\xb$ with $\yb$ assumed to be the origin.

A possibly more interesting argument in favor of studying two-point functions (especially in nonlocal settings) is that starting from a one-point function $\vb(\zbm)$ (usually provided), one can easily generate a two-point function $\ub(\xb,\yb)$ to which the decomposition results would apply, by simply taking
\[
\ub(\xb, \yb):=\frac{1}{2}\left(\vb(\xb)+\vb(\yb)\right),
\]
\[
\ub(\xb, \yb):=\vb(\xb-\yb),
\]
or
\[
\ub(\xb, \yb):=\vb(\xb) \psi(\xb-\yb),
\]
where $\psi$ is a given scalar function that possibly measures the interaction between $\xb$ and $\yb$. The choices here are infinite, however, the context of the application, physical or mathematical considerations could be invoked to choose a suitable candidate. In fact, for a single given one-point function, one may desire to obtain a set of decompositions for different two-point functions generated by the initial one, each of these two-point functions capturing a different type of interaction or behavior.

%%%%%%%%%%%%%%%%%%%%%%%%%%%%%%%%%%%%%%%%%%%%%%%%%%%%%%%%%%%%%%%%%%%%%%
%%%%%%%%%%%%%%%%%%%%%%%%%%%%%%%%%%%%%%%%%%%%%%%%%%%%%%%%%%%%%%%%%%%%%%
\section{Acknowledgements.}

The authors would like to thank the Michigan Center for Applied and Interdisciplinary Mathematics (MCAIM) at University of Michigan and AWM for supporting the WIMM workshop \cite{WIMM} where this work originated.

M. D'Elia was supported by Sandia National Laboratories (SNL), SNL is a multimission laboratory managed and operated by National Technology and Engineering Solutions of Sandia, LLC., a wholly owned subsidiary of Honeywell International, Inc., for the U.S. Department of Energys  National  Nuclear  Security  Administration  contract  number  DE-NA-0003525. This paper describes objective technical results and analysis. Any subjective views or opinions that might be expressed in the paper do not necessarily  represent the views of the U.S. Department of Energy or the United States Government. SAND Number: SAND2019-9727 O. X. Li was supported by NSF-DMS 1720245 and a UNC Charlotte faculty research grant. P. Radu was supported by NSF-DMS 1716790. Y. Yu was supported by NSF-DMS 1620434 and a Lehigh faculty research grant.


\begin{thebibliography}{10}

\bibitem{WIMM}
{WIMM} workshop, university of michigan, 2018.
\newblock \newline
  https://mcaim.math.lsa.umich.edu/events/women-in-mathematics-of-materials-wo%
rkshop/.

\bibitem{Mengesha2010}
B.~Aksoylu and T.~Mengesha.
\newblock Results on non-local boundary-value problems.
\newblock {\em Numerical Functional Analysis and Optimization}, 31:1301--1317,
  2010.

\bibitem{bhatia2013helmholtz}
Harsh Bhatia, Gregory Norgard, Valerio Pascucci, and Peer-Timo Bremer.
\newblock The {H}elmholtz-{H}odge decomposition: a survey.
\newblock {\em IEEE Transactions on visualization and computer graphics},
  19(8):1386--1404, 2013.

\bibitem{H2}
Shui-Nee Chow, Wuchen Li, and Haomin Zhou.
\newblock Entropy dissipation of {F}okker-{P}lanck equations on graphs.
\newblock {\em arXiv preprint arXiv:1701.04841}, 2017.

\bibitem{delia_fractional}
M.~D'Elia and M.~Gunzburger.
\newblock The fractional laplacian operator on bounded domains as a special
  case of the nonlocal diffusion operator.
\newblock {\em Computers and Mathematics with applications}, 66:1245--1260,
  2013.

\bibitem{sirev}
Q.~Du, M.~Gunzburger, R.~Lehoucq, and K.~Zhou.
\newblock Analysis and approximation of nonlocal diffusion problems with volume
  constraints.
\newblock {\em SIAM Review}, 54(4):667--696, 2012.

\bibitem{Du2013c}
Q.~Du, M.~Gunzburger, R.~Lehoucq, and K.~Zhou.
\newblock Analysis of the volume-constrained peridynamic navier equation of
  linear elasticity.
\newblock {\em Journal of Elasticity}, 113:193--217, 2013.

\bibitem{nlvc}
Q.~Du, M.~Gunzburger, R.~B. Lehoucq, and K.~Zhou.
\newblock A nonlocal vector calculus, nonlocal volume--constrained problems,
  and nonlocal balance laws.
\newblock {\em Mathematical Models and Methods in Applied Sciences},
  23(03):493--540, 2013.

\bibitem{Du2011b}
Q.~Du and K.~Zhou.
\newblock Mathematical analysis for the peridynamic non-local continuum theory.
\newblock {\em ESAIM: Mathematical Modelling and Numerical Analysis},
  45:217--234, 2011.

\bibitem{Emmrich2007b}
E.~Emmrich and O.~Weckner.
\newblock On the well-posedness of the linear peridynamic model and its
  convergence towards the navier equation of linear elasticity.
\newblock {\em Communications in Mathematical Sciences}, 5:851--864, 2007.

\bibitem{Emmrich2007a}
E.~Emmrich and O.~Weckner.
\newblock The peridynamic equation and its spatial discretization.
\newblock {\em Mathematical Modelling and Analysis}, 12:17--27, 2007.

\bibitem{fos-rad}
Mikil Foss and Petronela Radu.
\newblock Differentiability and integrability properties for solutions to
  nonlocal equations.
\newblock In {\em New Trends in Differential Equations, Control Theory and
  Optimization: Proceedings of the 8th Congress of Romanian Mathematicians},
  pages 105--119. World Scientific, 2016.

\bibitem{foss2019}
Mikil~D. Foss and Petronela Radu.
\newblock {\em Bridging Local and Nonlocal Models: Convergence and Regularity},
  pages 1243--1263.
\newblock Springer International Publishing, Cham, 2019.

\bibitem{FRW}
Mikil~D. Foss, Petronela Radu, and Cory Wright.
\newblock Existence and regularity of minimizers for nonlocal energy
  functionals.
\newblock {\em Differential Integral Equations}, 31(11/12):807--832, 11 2018.

\bibitem{Lehoucq-Gunzburger}
M.~Gunzburger and R.~B. Lehoucq.
\newblock A nonlocal vector calculus with application to nonlocal boundary
  value problems.
\newblock {\em Multiscale Modeling and Simulation}, 8:1581--1598, 2010.

\bibitem{Du2019}
H.~Lee and Q.~Du.
\newblock Nonlocal gradient operators with a nonspherical interaction
  neighborhood and their applications.
\newblock Technical Report arXiv:1903.06025, ArXiV, 2019.

\bibitem{gspl}
M.L.~Parks M.D.~Gunzburger, P.~Seleson and R.B. Lehoucq.
\newblock Peridynamics as an upscaling of molecular dynamics.
\newblock {\em Multiscale Modeling and Simulation}, (8):204--227, 2009.

\bibitem{DuMengesha}
Tadele Mengesha and Qiang Du.
\newblock The bond-based peridynamic system with {D}irichlet-type volume
  constraint.
\newblock {\em Proc. Roy. Soc. Edinburgh Sect. A}, 144(1):161--186, 2014.

\bibitem{RTY}
Petronela Radu, Daniel Toundykov, and Jeremy Trageser.
\newblock A nonlocal biharmonic operator and its connection with the classical
  analogue.
\newblock {\em Archive for Rational Mechanics and Analysis}, 223(2):845--880,
  2017.

\bibitem{radu_wells2019}
Petronela Radu and Kelsey Wells.
\newblock A doubly nonlocal {Laplace} operator and its connection to the
  classical {L}aplacian.
\newblock {\em J. Integral Equations Applications}, 2019.
\newblock Advance publication.

\bibitem{Silling2000}
S.A. Silling.
\newblock Reformulation of elasticity theory for discontinuities and long-range
  forces.
\newblock {\em Journal of the Mechanics and Physics of Solids}, 48:175--209,
  2000.

\bibitem{silling2010peridynamic}
SA~Silling and RB~Lehoucq.
\newblock Peridynamic theory of solid mechanics.
\newblock {\em Advances in applied mechanics}, 44:73--168, 2010.

\bibitem{Du2010c}
K.~Zhou and Q.~Du.
\newblock Mathematical and numerical analysis of linear peridynamic models with
  non-local boundary.
\newblock {\em SIAM Journal on Numerical Analysis}, 48:1759--1780, 2010.

\end{thebibliography}
\end{document}